\documentclass[12pt,a4paper,twoside]{article}
\usepackage{a4wide,amsmath,amsbsy,amsfonts,amssymb,stmaryrd}
\newcommand\ZZ{{\hat{\mathbb Z}}}
\newcommand\Z{{\mathbb Z}}

\newcommand\Q{{\mathbb Q}}

\newcommand\C{{\mathbb C}}

\newcommand\ra{\rightarrow}

\newcommand\ilim{\lim\limits_{\longleftarrow}\,}

\newcommand\SL{\mathrm{SL}}

\newcommand\SP{\mathrm{Sp}}

\newcommand\aut{\mathrm{Aut}}
\newcommand\out{\mathrm{Out}}
\newcommand\inn{\mathrm{inn}}
\newcommand\spec{\mathrm{Spec}}
\renewcommand\dim{\mathrm{dim}}

\newcommand\hookra{\hookrightarrow}
\newcommand\tura{\twoheadrightarrow}
\newcommand\da{\downarrow}

\newcommand\dd{\partial}
\renewcommand{\hom}{\mathrm{Hom}}

\renewcommand{\Im}{\mathrm{Im}} 
\newcommand\sr{\stackrel}
\newcommand\st{\scriptstyle}
\newcommand\sst{\scriptscriptstyle}
\newcommand\hGG{\hat{\GG}}
\newcommand\hPi{\hat{\Pi}}

\newcommand\hpi{\hat{\pi}}
\newcommand\cpi{\check{\pi}}

\newcommand\ssm{\smallsetminus}
\newcommand\ol{\overline}
\newcommand\ccM{\overline{\cal M}}
\newcommand\cM{{\cal M}}
\newcommand\cA{{\cal A}}
\newcommand\cB{{\cal B}}

\newcommand\cK{{\cal K}}
\newcommand\cG{{\cal G}}

\newcommand\cU{{\cal U}}
\newcommand\cT{{\cal T}}
\newcommand\cJ{{\cal J}}
\newcommand\GG{\Gamma}
\newcommand\ld{\lambda}

\newcommand\wh{\widehat}
\newcommand\wt{\widetilde}

\newcommand\wM{\wt{\cM}}
\newcommand\hM{\wh{\cM}}

\newcommand\hd{\hat{\Delta}}

\newcommand\td{\tilde}
\newcommand\sg{\sigma}
\newcommand\Sg{\Sigma}
\newcommand\gm{\gamma}

\newtheorem{theorem}{Theorem}[section]
\newtheorem{corollary}[theorem]{Corollary}
\newtheorem{proposition}[theorem]{Proposition}
\newtheorem{lemma}[theorem]{Lemma}

\newtheorem{remark}[theorem]{Remark}

\newenvironment{prf}[1]{\trivlist
\item[\hskip
\labelsep{\it #1.\hspace*{.3em}}]}{~\hspace{\fill}~$\square$\endtrivlist}
\newenvironment{proof}{\begin{prf}{Proof}}{\end{prf}}

\begin{document}

\title{Fundamental groups of moduli stacks \\ of stable curves of compact type}
\author{Marco Boggi}\maketitle

\begin{abstract}
Let $\wM_{g,n}$, for $2g-2+n>0$, be the moduli stack of $n$-pointed, genus $g$, stable complex curves
of compact type. Various characterizations and properties are obtained of both the algebraic and
topological fundamental groups of the stack $\wM_{g,n}$. 

Let $\GG_{g,n}$, for $2g-2+n>0$, be the Teichm{\"u}ller group associated with a compact Riemann surface of
genus $g$ with $n$ points removed $S_{g,n}$, i.e. the group of homotopy classes of diffeomorphisms of
$S_{g,n}$  which preserve the orientation of $S_{g,n}$ and a given order of its punctures. Let $\cK_{g,n}$ be
the normal subgroup of $\GG_{g,n}$ generated by Dehn twists along separating circles on $S_{g,n}$.
As a first application of the above theory, a characterization of $\cK_{g,n}$ is given for all $n\geq 0$ (for
$n=0,1$, this was done by Johnson in \cite{J3}). We define the {\it Torelli} group $\cT_{g,n}$ to be the
kernel of the natural representation $\GG_{g,n}\ra\SP_{2g}(\Z)$. The abelianization of the Torelli group
$\cT_{g,n}$ is determined for all $g\geq 1$ and $n\geq 1$, thus completing classical results by 
Johnson~\cite{J4} and Mess~\cite{Mess} (for a different definition of the Torelli group, this 
was done by van den Berg in \cite{vdb}, who provides also a new proof of Johnson's and Mess' results).
\newline

\end{abstract}

\section{Introduction}\label{intro}
Let $\cM_{g,n}$ be the moduli stack of smooth, $n$-pointed, genus $g$, complex curves and let
$\ccM_{g,n}$ be its Deligne-Mumford compactification, i.e. the moduli stack of stable, $n$-pointed, 
genus $g$, complex curves. They are actually Deligne-Mumford stacks (briefly D--M stacks).
A central role in the study of the space $\cM_{g,n}$ (and then of $\ccM_{g,n}$ as well) is played by the
Teichm{\"u}ller group $\GG_{g,n}$. Classically, it is defined, for a given compact Riemann surface of genus $g$
with $n$ points removed $S_{g,n}$, as the {\it mapping class group} of $S_{g,n}$, i.e. the group of homotopy
classes of diffeomorphisms of $S_{g,n}$ which preserve the orientation of $S_{g,n}$ and a given order of its
punctures. From our point of view, the most significant characterization of $\GG_{g,n}$ is as the
fundamental group of the topological stack underlying $\cM_{g,n}$. Since, topologically, the space
$\cM_{g,n}$ is an Eilenberg-MacLane space of type $(\pi,1)$, in principle, all topological properties of
$\cM_{g,n}$ can be derived from the study of $\GG_{g,n}$. According to the anabelian philosophy of
Grothendieck (see \cite{B-L}), much more is true (see also the rigidity results in \cite{B2}).
This undoubtedly accounts for the extraordinary richness and difficulty of Teichm{\"u}ller theory.

At the other extreme lies the case of the D--M stack $\ccM_{g,n}$, whose fundamental group is trivial
(see Proposition~1.1 in \cite{B-P}).

The purpose of this paper is a sistematic study of the topological and algebraic fundamental group of a
partial compactification of $\cM_{g,n}$, contained as an open sub-stack in the Deligne-Mumford
compactification $\ccM_{g,n}$, which is, in some sense, intermediate between the two. Let
$\wM_{g,n}$ be the moduli stack of $n$-pointed, genus $g$ stable complex curves of compact type. In 
other words, the stack $\wM_{g,n}$ parametrizes stable curves which have a compact generalized 
jacobian or, equivalently, whose dual graph is a tree. 

The topological and algebraic fundamental groups of $\wM_{g,n}$ are characterized in various ways. Let
$\cK_{g,n}$ be the normal subgroup of $\GG_{g,n}$ generated by Dehn twists along separating circles on
$S_{g,n}$. The first (and almost trivial) characterization of $\pi_1(\wM_{g,n})$ is as the quotient of the
mapping class group $\GG_{g,n}$ by $\cK_{g,n}$. Eventually, it is proved that $\pi_1(\wM_{g,n})$ is a linear
group and, more precisely, an extension of the symplectic group $\SP_{2g}(\Z)$ by a free abelian group  of
rank ${2g\choose 3}+{\st 2g(n-1)}$. This generalizes to all $g\geq 1$ and $n\geq 1$ previous results by
Morita, Hain and Looijenga (see \cite{M1}  and \cite{H-L}). In particular, it turns out to be a much simpler
object than the Teichm{\"u}ller group itself. It is important to stress that, in contrast with the usual
situation, sometimes, the characterizations of the topological fundamental group are derived from those of
the algebraic one (this is the case for Theorem~\ref{johnson's kernel}, which is the key result of the paper).

The above theory turns out to be very useful because, even though the space $\wM_{g,n}$
in general is not an Eilenberg-MacLane space of type $(\pi,1)$ (see Proposition~2.4 in \cite{Mondello}), its
fundamental group still contains a lot of informations on the moduli stack $\wM_{g,n}$ and on
its open substack $\cM_{g,n}$ as well. This will be apparent from the applications
given in  \cite{B2} to the monodromy of families of curves of compact type. 
As to $\cM_{g,n}$, let us just mention Proposition~\ref{morita} and Remark~\ref{kawazumi}. 

However, the most important applications, given in the present paper, of the above results are to some
questions in classical Teichm{\"u}ller theory. For $n=0$ and $g\geq 2$, Johnson charecterized $\cK_{g,n}$ as
the kernel of the natural representation:
$$\GG_g\ra\out(\pi_1(S_g)/\pi_1(S_g)^{[3]}),$$
where $\pi_1(S_g)^{[3]}$ is the third term of the descending central series of the fundamental group of
$S_g$. From the point of view of Hodge theory, the {\it weight filtration} $W^k\pi_1(S_{g,n})$, for $k\geq 0$,
(see Section~\ref{level} for the definition) is the natural extension of the descending central series to
the $n$-punctured case. In Theorem~\ref{johnson's kernel}, it is proved that, for $2g-2+n>0$, the group
$\cK_{g,n}$ is the kernel of the natural representation:
$$\GG_{g,n}\ra\out(\pi_1(S_{g,n})/W^3\pi_1(S_{g,n})).$$

Let $\cT_{g,n}$ be the Torelli subgroup of $\GG_{g,n}$. It is commonly defined as the kernel of the natural representation $\GG_{g,n}\ra\SP_{2g}(\Z)$, i.e. as the subgroup of mapping classes of $\GG_{g,n}$ acting 
trivially on the homology of the compact surface $S_g$. Even though, for $n\geq 2$, there is no complete
agreement on this definition, it is certainly the most natural from the point of view of algebraic geometry,
i.e. Hodge theory, and we stick to it. In a series of papers (see \cite{J1}, \cite{J2}, \cite{J3}, \cite{J4}), Johnson 
computed the abelianization of $\cT_{g,n}$ for $g\geq 3$ and $n=0,1$, while Mess computed the abelianization 
of $\cT_2$ (see \cite{Mess}). Then in \cite{H-L}, Hain and Looijenga implemented their results in order to
compute $H_1(\cT_{g,n},\Q)$ for all $g\geq 1$ and $n\geq 1$. Thanks to the above characterization of the
group $\cK_{g,n}$, I complete the above results, determining the abelianization of the Torelli group
$\cT_{g,n}$ for all $g\geq 1$ and $n\geq 1$ (see Theorem~\ref{abelianization}).
 
Let us in the end mention a result which, in the light of \cite{B2}, can also be considered one 
of the most important of the paper. For $g_1$ and $g_2$ non-negative integers, such that $g_1+g_2=g$, and a
partition of the set of the $n$ marking labels in two subsets of cardinality $n_1$ and $n_2$, Knudsen
defined a clutching morphism $\ccM_{g_1,n_1+1}\times\ccM_{g_2,n_2+1}\ra\ccM_{g,n}$ 
(see \cite{Knudsen}), which, on the moduli space of stable curves of
compact type, restricts to a morphism $\wM_{g_1,n_1+1}\times\wM_{g_2,n_2+1}\ra\wM_{g,n}$. In
Theorem~\ref{stratum}, it is proved that such morphism induces a monomorphism between the respective
topological fundamental groups (for the algebraic fundamental groups the situation is a little bit trickier,
due to the failure of the congruence subgroup property for $\SL_2(\Z)$). In particular, the operad structure
of the category, with objects the stacks $\wM_{g,n}$ and with arrows the associated clutching morphisms,
is faithfully reflected by the operad structure of the category of their topological fundamental groups with
the induced homomorphisms.

\section{Level structures over moduli of curves}\label{level}
The purpose of this section is, basically, to provide the notations to be used in the rest of the paper.
For a more complete treatment of level structures and Teichm{\"u}ller theory, we refer the reader, for
instance, to \S 1 of \cite{B}.

Let $\ccM_{g,n}$ (with $2g-2+n>0$) be the stack of $n$-pointed, genus $g$, stable algebraic
curves over $\C$ (see \cite{DM}).  It is a regular connected proper Deligne--Mumford stack (briefly {\it
D--M stack}) over $\C$ of dimension $3g-3+n$, and it contains, as an open substack, the stack $\cM_{g,n}$ 
of $n$-pointed, genus $g$, smooth algebraic curves over $\C$. The stack $\ccM_{g,n}$ is simply connected
(see Proposition~1.1 in \cite{B-P}). On the contrary, the stack $\cM_{g,n}$ has plenty of non-trivial covers
which we are briefly going to introduce in this section. Its universal cover, the Teichm{\"u}ller space
$T_{g,n}$, is a contractible complex manifold. There is a natural way to define homotopy groups for
topological D--M stacks (see \cite{N1} and \cite{N2}). Then, the fundamental group $\pi_1(\cM_{g,n},a)$,
called the Teichm{\"u}ller group and denoted by $\GG_{g,n}$,  is isomorphic to the deck 
transformations' group of the cover $T_{g,n}\ra\cM_{g,n}$. 

A level structure $\cM^\ld$ is a finite, connected, Galois, {\'e}tale cover of the stack 
$\cM_{g,n}$ (here, an {\'e}tale cover is an {\'e}tale, surjective, representable morphism of
algebraic stacks), therefore it is also represented by a regular Deligne--Mumford stack
$\cM^\ld$.  The {\it level} associated to $\cM^\ld$ is the finite index normal subgroup
$\GG^\ld:=\pi_1(\cM^\ld, a')$ of the Teichm{\"u}ller group $\GG_{g,n}$.
 
A level structure $\cM^{\ld'}$ {\it dominates} $\cM^\ld$, if there is a natural {\'e}tale
morphism $\cM^{\ld'}\ra\cM^{\ld}$ or, equivalently, $\GG^{\ld'}\leq\GG^\ld$. To mark the fact that
$\cM^\ld$ is a level structure over $\cM_{g,n}$, we will often denote it by $\cM_{g,n}^\ld$.
  
The morphism $p:\cM_{g,n+1}{\ra}\cM_{g,n}$ (forgetting the last label) is
naturally isomorphic to the universal curve over $\cM_{g,n}$. One can then identify the fiber
$p^{-1}(a)$ with an $n$-punctured, genus $g$ curve $C\ssm\{P_1,\ldots,P_n\}$, where
$(C,P_1,\ldots,P_n)$ is a curve in the class determined by $a\in\ccM_{g,n}$. Denote by $S_{g,n}$ the
Riemann surface underlying $C\ssm\{P_1,\ldots,P_n\}$ and fix a point $\tilde{a}\in S_{g,n}$. Since
$p$ is a Serre fibration and $\pi_2(\cM_{g,n})=\pi_2(T_{g,n})=0$, there is a short exact sequence on
fundamental groups
$$1\ra\pi_1(S_{g,n},\tilde{a})\ra\pi_1(\cM_{g,n+1},\tilde{a})\ra\pi_1(\cM_{g,n},a)\ra 1.$$
By a standard argument this defines a monodromy representation:
$$\rho:\pi_1(\cM_{g,n},a)\ra \mbox{Out}(\pi_1(S_{g,n},\tilde{a})),$$
called the {\it universal monodromy representation}.
From the above fibration, it 
follows also that there is a natural representation of $\GG_{g,n}$ in the group of homotopy (or
equivalently isotopy) classes of self-homeomorphism of the Riemann surface $S_{g,n}$. Let us
denote by $\hom^+(S_{g,n})$ the subgroup of orientation preserving self-homeomorphism of $S_{g,n}$
and by $\hom^0(S_{g,n})$ the subgroup consisting of homeomorphism homotopic (isotopic) to the
identity. By a classical result in Teichm{\"u}ller theory,  the representation associated to the
universal curve actually induces an isomorphism of the Teichm{\"u}ller group with the so called {\it
mapping class group} of $S_{g,n}$: 
$$\GG_{g,n}\cong\hom^+(S_{g,n})/\hom^0(S_{g,n}).$$  

Let us denote by $\Pi_{g,n}$ the fundamental group of $S_{g,n}$ based in
$\tilde{a}$ and give $\Pi_{g,n}$ the standard presentation:
$$ \Pi_{g,n}=<\alpha_1, \dots \alpha_g,\beta_1,\dots,\beta_g, u_1,\dots,u_n|\;
\prod_{i=1}^g[\alpha_i,\beta_i] \cdot u_n\cdots u_1>,$$ 
where $u_i$, for $i=1,\dots,n$, is a simple loop around the  puncture $P_i$. For $n\ge 1$, let
$A(g,n)$ be the group of automorphisms of $\Pi_{g,n}$ which  fix the conjugacy classes of
all $u_i$. For $n=0$,  let instead $A(g,0)$ be the image of
$A(g,1)$ in the automorphism group of $\Pi_g:=\Pi_{g,0}$. Finally, let $I(g,n)$ be the group of 
inner  automorphisms of $\Pi_{g,n}$.  With these notations, the representation $\rho$ is
faithful and gives a natural isomorphism $\GG_{g,n}\cong A(g,n)/I(g,n)$.

The most natural way to define levels is provided by the above isomorphism. In general,
for a subgroup $\Pi^\ld\le\Pi_{g,n}$, which is invariant under $A(g,n)$ (in such
case we simply say that $\Pi^\ld$ is {\it invariant}), it is defined a representation:
$$\rho_\ld:\GG_{g,n}\ra \mbox{Out}(\Pi_{g,n}/\Pi^\ld),$$
whose kernel we denote by $\GG^\ld$. When $\Pi^\ld$ has finite index in $\Pi_{g,n}$, then
$\GG^\ld$ has finite index in $\GG_{g,n}$ and is called the {\it geometric level} associated to
$\Pi^\ld$. The corresponding level structure is denoted by $\cM_{g,n}^{\ld}$. 

A class of finite index invariant subgroups of the group $\Pi_{g,n}$ one can consider, in order to
define geometric level structures, is that obtained from the descending central series, twisting by $l$-th
powers. The descending central series is defined by $\Pi^{[1]}:=\Pi_{g,n}$ and
$\Pi^{[k]}:=[\Pi^{[k-1]},\Pi]$. Let then $\Pi^l$ be the invariant subgroup of $\Pi_{g,n}$
spanned by $l$-th powers and define 
$$\Pi^{[k],l}:= \Pi^{[k]}\cdot\Pi^l.$$
We denote the corresponding levels and level structures by $\GG^{[k],l}$ and $\cM^{[k],l}$. 

Even though algebraically, the above filtration is very natural, from a geometric
point of view the so called {\it weight filtration} is more significative. This filtration originates in
Hodge theory and is defined as follows. Let $N$ be the kernel of the natural morphism
$\Pi_{g,n}\ra\Pi_{g,0}$ (filling in the punctures) and define
\[\begin{array}{lll}
W^1 \Pi& := &\Pi_{g,n},\\
W^2 \Pi& := &N\cdot \Pi^{[2]},\\
W^{k+1} \Pi& := &[\Pi_{g,n},W^{k} \Pi]\cdot [N, W^{k-1}\Pi].
\end{array}\]
As for the descending central series, one has $[W^s \Pi,W^t \Pi]\leq W^{s+t} \Pi$.
The descending central series and the weight filtration are cofinal to each other (for $n=0$, they
coincide). In fact $W^{2k-1}\Pi\leq\Pi^{[k]}\leq W^k \Pi$. As above, let us define
$$W^{k,l} \Pi=W^{k} \Pi\cdot\Pi^l.$$
The corresponding representation is denoted by $\rho_{w(k,l)}$
and the corresponding levels and level structures are denoted by $\GG^{w(k,l)}$ and $\cM^{w(k,l)}$
respectively. The kernel of the representation
$\rho_{W^k}:\GG_{g,n}\ra\mbox{Out}(\Pi_{g,n}/W^k\Pi)$ will be instead denoted by $W^k\GG$.

Of particular interest are the levels defined by the kernels of the representations:
$$\rho_{(m)}:\GG_{g,n}\ra \mbox{Sp}(H_1(S_g,\Z/m)),\;\;\mbox{ for }m\geq 2.$$
They are denoted by $\GG(m)$ and called {\it abelian levels of order $m$}. The corresponding
level structures are then denoted by $\cM^{(m)}$.

The kernel  of the representation $\GG_{g,n}\ra\mbox{Sp}_{2g}(\Z)$ is denoted by $\cT_{g,n}$
and called the {\it Torelli subgroup of $\GG_{g,n}$.} Note that $\GG(m)=\GG^{w(2,m)}$ and
$\cT_{g,n}=W^2\GG$. 

The usual way to compactify a level structure $\cM^{\ld}$
over $\cM_{g,n}$ is to take the normalization of $\ccM_{g,n}$ in the function field of $\cM^{\ld}$.
A more functorial definition can be given in the category of log regular schemes. Let $\dd$ be the
logarithmic structure on $\ccM_{g,n}$ associated to the normal crossing divisor
$\dd\cM:=\ccM_{g,n} \ssm\cM_{g,n}$. A level structure over $(\ccM_{g,n},\dd)$ is a finite,
connected, log {\'e}tale cover 
$$(\ccM^\ld,\dd^\ld)\ra(\ccM_{g,n},\dd).$$
The log purity Theorem (see \cite{Moch}) implies that any level structure $\cM^\ld$ over
$\cM_{g,n}$ admits a canonical compactification to a level structure $(\ccM^\ld,\dd^\ld)$ over
$(\ccM_{g,n},\dd)$, where $\ccM^\ld$ is the normalization of $\ccM_{g,n}$ in the function field of
$\cM^\ld$ and $\dd^\ld$ the logarithmic structure associated to the normal crossing divisor $\dd
\cM^\ld:=\ccM^\ld\ssm\cM^\ld$ (the {\it Deligne-Mumford boundary} of
$\ccM^{\ld}$). On the other hand, it is also clear that any level structure over $(\ccM_{g,n},\dd)$ 
can be realized in this way. So, forgetting the logarithmic structure, one is  back to the 
previous definition. The following result is classical:
 
\begin{proposition}\label{Deligne}If a level $\GG^\lambda$ is contained in an abelian level
of order $m$, for some $m\ge 3$, then the level structure $\ccM^{\ld}$ is
represented by a projective variety.
\end{proposition}
   
There is a very explicit and elementary method to describe the compactifications $\ccM^\ld$, locally in the
analytic topology.  A neighborhood of a point $a\in\cM^{\ld}$ is just the base of the local universal
deformation of the fibre $C_a$ in $a$ of the universal family ${\cal C}^\ld\ra\cM^\ld$. Let us see how a
neighborhood of $a\in\dd\cM^\ld$ can be described.      

Let $\cB\ra\ccM_{g,n}$ be an analytic neighborhood of the image $y$ of $a$ in $\ccM_{g,n}$ such
that:
\begin{itemize} 
\item local coordinates $z_1,...,z_{3g-3+n}$ embeds $\cB$ in $\C^{3g-3+n}$ as an open ball;
\item $C:=\pi^{-1}(y)$ is the most degenerate curve in the pull-back ${\cal C}\stackrel{\pi}{\ra}\cB$
of the universal family over $\cB$; 
\item  an {\'e}tale groupoid representing $\ccM_{g,n}$ trivializes over $\cB$ to 
$\mbox{Aut}(C)\times\cB\rightrightarrows\cB$. 
\end{itemize}
Let $\{Q_1,\dots,Q_s\}$ be the set of singular points of $C$ and let
$z_i$, for $i=1,\dots ,s$, parametrize curves where the singularity $Q_i$
subsists. The discriminant locus $\dd\cB\subset\cB$ of $\pi$
has then equation $z_1\cdots z_s=0$. Let $U=\cB\ssm\dd\cB$. The natural morphism $U\ra\cM_{g,n}$ 
induces a homomorphism of fundamental groups:
$$\phi_\ld:\pi_1(U,a)\ra\pi_1(\cM_{g,n},a).$$ 
If we denote as well by $\GG^\ld$ the subgroup of $\pi_1(\cM_{g,n},a)$ determined by the {\'e}tale cover
$\cM^\ld\ra\cM_{g,n}$, a connected component $U^{\ld}$ of $U\times_{\ccM_{g,n}}\ccM^{\ld}$ is then
determined by the subgroup $\phi_\ld^{-1}(\GG^\ld)$ of the abelian group $\pi_1(U,a)$.

Let us make the above description more explicit.
Let $s_1,\dots, s_n$ be the sections of the universal family over
$\ccM_{g,n}$ and define ${\cal C}|_U:={\cal C}\ssm\bigcup_{i=1}^n s_i(U)$.
Let us fix a homeomorphism between the fiber over $a$ of the morphism ${\cal C}|_U\ra U$ and the
marking Riemann surface $S_{g,n}$. In more detail, the local monodromy representation
$$\rho_U:\pi_1(U,a)\ra \mbox{Out}(\pi_1(S_{g,n},\tilde{a})),$$  
associated to the family ${\cal C}|_U\ra U$, is defined as the composition of the natural morphism
$\pi_1(U)\ra\pi_1(\cM_{g,n})$ with the universal monodromy representation and can be explicitly 
described as follows.
    
Since $U$ is homotopic to the $s$-dimensional torus $(S^1)^s$, the fundamental
group $\pi_1(U,a)$ is abelian and freely generated by simple loops $\gamma_i$
around the divisors $z_i$, for $i=1,\dots,s$. Such loops can be lifted to disjoint loops
$\tilde{\gamma}_i$, for $i=1,\ldots,s$, in $S_{g,n}$, whose isotopy classes are uniquely
determined and which become isotrivial after specialization to $C$. One can prove that 
$\gamma_i$ is mapped by $\rho_U$ exactly in the element of Out$(\pi_1(S_{g,n},\tilde{a}))$
determined by the Dehn twist $\tau_{\tilde{\gamma}_i}$ along $\tilde{\gamma}_i$, for
$i=1,\dots,s$. In particular, the representation $\rho_U$ is faithful. 

Let $E_{\Sigma(C)}$ be the free abelian group spanned by the edges of the dual graph $\Sigma(C)$ of
the stable curve $C$. The edges of the dual graph correspond to  isotopy classes of circles $\gm_e$ in 
$S_{g,n}$ which become isotrivial specializing to $C$. The group $E_{\Sigma(C)}$ can then be naturally
identified with the free abelian group spanned in $\GG_{g,n}$ by the set of Dehn twists $\{\tau_{\gm_e}\}$.
On the other hand, the fundamental group $\pi_1(U,a)$ is as well naturally isomorphic to $E_{\Sigma(C)}$. So,
for a   level $\GG^\ld$ the subgroup $\phi_\ld^{-1}(\GG^\ld)$ of the fundamental group $\pi_1(U,a)$ is then
canonically identified with $E_{\Sigma(C)}\cap\GG^\ld$.

An almost complete description of local monodromy coefficients for the geometric levels $\GG^{w(k,l)}$ 
is given in Theorem~3.1.3, \cite{P-dJ}, Proposition~2.8, \cite{B-P}, and Theorem~3.3.3, \cite{Pikaart}. 
Let us collect all their results in a single statement. 

Let $N_{\Sigma(C)}$ and $S_{\Sigma(C)}$ be respectively the subgroups of 
$E_{\Sigma(C)}$ spanned by edges corresponding to  non-separating circles and by edges corresponding to
separating circles. Let then $S^1_{\Sg(C)}$ be the subgroup of $S_{\Sigma(C)}$ spanned by edges
corresponding to separating circles bounding an unpunctured genus $1$ surface. Let instead $P_{\Sigma(C)}$
be the subgroup spanned  by elements of the form $e_1-e_2$, where $\{e_1,e_2\}$ corresponds to a cut 
pair on $C$. Eventually,  denote by $P_{\Sigma(C)}^{un}$ the subgroup of $P_{\Sigma(C)}$ spanned by
elements corresponding to cut pairs bounding a surface without punctures. Let, for $m,s$ positive integers,
$m_s:=m/\gcd(m,s)$.

\begin{theorem}\label{monodromy coefficients}With the above notations, 
the kernel of $\rho_U^{w(k,l)}$, where $U$ is a neighborhood of $[C]\in\ccM_{g,n}$ as above,
is given by:
\[\begin{array}{lll}
1.&\mbox{If $k=2$:} &m\, N_{\Sg(C)}+P_{\Sigma(C)}+ S_{\Sigma(C)}. \\
2.&\mbox{If $k=3$:}&m\,N_{\Sg(C)}+m_2\,P_{\Sg(C)}^{un}+S_{\Sg(C)}.\\ 
3.&\mbox{If $k=4$ and $m$ is odd or $4|m$: \,\,}&
m\,N_{\Sg(C)}+m_2\,S_{\Sg(C)}+m_6\,S^1_{\Sg(C)}.\\
&\mbox{If $k=4$ and $2||m$:}&m\,N_{\Sg(C)}+m_2\,P_{\Sg(C)}+
m_2\,S_{\Sg(C)}+m_6\,S^1_{\Sg(C)}.\\ 
4.&\mbox{If $k\ge 4$ and gcd$(m,6)=1$:\,\,}&m\,E_{\Sg(C)}.
\end{array} \]
\end{theorem}

\section{Moduli spaces of stable curves of compact type}\label{compact type}
In this section, we begin the study of the moduli space $\wM_{g,n}$ of stable $n$-pointed, genus
$g$, complex curves of compact type. We will give various characterization both of its topological
fundamental group $\pi_1(\wM_{g,n},a)$ and of its algebraic fundamental group $\hpi_1(\wM_{g,n},a)$
(in the sequel, we will often omit any mention of base points). As usual, the profinite completion
of a given group $G$ is denoted by $\hat{G}$ (this motivates the previous notation for the algebraic
fundamental group of $\wM_{g,n}$).

Since $\dd\wM_{g,n}:=\wM_{g,n}\ssm\cM_{g,n}$ is a normal crossing divisor, the embedding
$\cM_{g,n}\hookra\wM_{g,n}$ induces an epimorphism $\pi_1(\cM_{g,n},a)\tura\pi_1(\wM_{g,n},a)$ on
fundamental groups, whose kernel is normally generated by small loops around the irreducible components
of $\dd\wM_{g,n}$. As we remarked in Section~\ref{level}, the isomorphism $\GG_{g,n}=\pi_1(\cM_{g,n},a)$
identifies such loops with Dehn twists along separating circles. Let us then define 
$$\cK_{g,n}:=<\tau_\alpha\in\GG_{g,n} |\, \alpha\mbox{ is a separating circle on }S_{g,n}>,$$
which is also called the {\it Johnson subgroup of} $\GG_{g,n}$.
We then get the first and the most immediate of the many characterizations we will give of
$\pi_1(\wM_{g,n},a)$:

\begin{proposition}\label{twists}Let $2g-2+n>0$. The fundamental group of $\wM_{g,n}$ fits in the short
exact sequence:
$$1\ra\cK_{g,n}\ra\GG_{g,n}\ra\pi_1(\wM_{g,n},a)\ra 1.$$
\end{proposition}

By Theorem~2 in \cite{Birman} (for the genus $2$ case) and Theorem~6 in \cite{J3} (for the genus $>2$ case),
we know that $\cK_g$ is the kernel, for $n=0$, of the natural representation
$$\rho_{W^3}: \GG_{g,n}\ra\out(\Pi_{g,n}/W^3\Pi).$$
In particular, the fundamental group of $\wM_g$ is naturally isomorphic to
$\Im\rho_{W^3}$. Passing to profinite completions, the above representation induces a continuous
representation:
$$\hat{\rho}_{W^3}: \hGG_{g,n}\ra\out(\hPi_g/W^3\hPi).$$
A natural guess then is that, for $n=0$, the algebraic fundamental group of $\wM_{g,n}$ is
isomorphic to $\Im\hat{\rho}_{W^3}$. Thanks to Corollary~3.11 in \cite{B-P} (basically), it is possible to
prove that this actually holds for all $n\geq 0$. In its turn, this will yield that the
fundamental group of $\wM_{g,n}$ is isomorphic to $\Im\rho_{W^3}$, for all $n\geq 0$.
Here, I will give more details of the proofs in \cite{B-P} and more complete results. Some
preliminary lemma's are needed.

\begin{lemma}\label{exact sequence} 
Let $f: {\cal C}\ra \cU$ be a flat, generically smooth curve over a connected 
analytic stack $\cU$.
Let ${\cal C}_u$ be the fibre over a point $u\in\cU$, and let also $\td{u}\in{\cal C}_u$.
There is then an exact sequence of fundamental groups:
$$\pi_1({\cal C}_u,\td{u})\ra\pi_1({\cal C},\td{u})\ra\pi_1(\cU,u)\ra 1.$$
\end{lemma} 
\begin{proof} It is clear that there is an epimorphism $f_*:\,\pi_1({\cal C})\ra
\pi_1(\cU)$, and that the kernel of $f_*$ is generated by all the loops 
contained in the fibres of $f$. One needs to show then that the loops can be
displaced horizontally from one fibre to the other. This can be done over
the open dense set in $\cU$ where $f$ is smooth, because there the family
is, topologically, locally trivial. But this suffices,
because in a flat family of curves there are specialization morphisms from
smooth to singular fibres, inducing epimorphisms on fundamental groups.

\end{proof}

Let us denote by $\beta_0:\ccM_{g-1,n+2}\ra\ccM_{g,n}$ the natural pinching morphism whose image is
the locus whose generic point parametrizes singular irreducible curves. An {\'e}tale cover $\wM^\ld$ of
$\wM_{g,n}$ compactifies to an {\'e}tale cover $\ccM^\ld\ra\ccM_{g,n}$ whose branch locus is contained
in the image of $\beta_0$. 

\begin{lemma}\label{comparison lemma}\begin{enumerate}
\item Let $\ccM^\ld\ra\ccM_{g,n}$ be a level structure whose branch locus is contained in the boundary
divisor $\beta_0(\ccM_{g-1,n+2})$. Denote by $p:\Gamma_{g,n}\ra\Gamma_{g,n-1}$ the
epimorphism induced filling in $P_n$ on $S_{g,n}$ and let $\GG^{p(\ld)}:=p(\Gamma^\ld)$. Then
$\pi_1(\ccM_{g,n}^\ld)=\pi_1(\ccM_{g,n-1}^{p(\ld)})$. 
\item Let $\GG^{\ld_1}\leq\GG^{\ld_2}$ be two levels whose associated level structures satisfy the 
hypothesis in the above item. If the natural morphism $\ccM_{g,n}^{\ld_1}\ra\ccM_{g,n}^{\ld_2}$ is {\'e}tale
and, with the same notations as above, moreover, $p(\Gamma^{\lambda_1})=p(\Gamma^{\lambda_2})$,
then it holds $\GG^{\ld_1}=\GG^{\ld_2}$.
\end{enumerate}
\end{lemma}

\begin{proof}Item (ii) follows directly from (i). So let us prove (i). We claim that the natural
morphism ${\cM}_{g,n}^\ld\stackrel{\phi}{\ra} {\cM}_{g,n-1}^{p(\ld)}$ is a fibration in smooth
curves. There is a factorization
\[\begin{array}{ccc}
 {\cM}_{g,n}^\ld &\stackrel{\tilde{\phi}}{\ra} &\cM_{g,n}\times_{\cM_{g,n-1}}\cM_{g,n-1}^{p(\ld)}\\
&\stackrel{\,\;\;\,\phi}{\searrow}&\da{\st \pi}\;\;\;\;\\ 
& &{\cM}_{g,n-1}^{p(\ld)},\end{array}\]
with $\tilde{\phi}$ {\'e}tale and $\pi$ smooth. Thus $\phi$ is also smooth.
In the Stein factorization of $\phi$
$${\cM}_{g,n}^\ld\stackrel{\phi'}{\ra}Y\stackrel{f}{\ra}{\cM}_{g,n-1}^{p(\ld)},$$
the morphism $\phi'$ has connected fibres and $f$ is finite. Therefore $f$ is {\'e}tale.
By definition of $\GG^{p(\ld)}$, the map induced on fundamental groups
by $\phi$ is surjective. So the same is true for the morphism $f$ which then is
an isomorphism. Hence $\phi=\phi'$.
In particular, since the morphism $\phi:\ccM_{g,n}^\ld\ra\ccM_{g,n-1}^{p(\ld)}$ is also log-smooth, it is a 
flat family of semi-stable curves. Let $S$ be a fibre of $\phi$ lying above a fibre of
$\ccM_{g,n}\ra\ccM_{g,n-1}$ which is a tree of $g$ elliptic curves.  By the assumptions made on 
the branch locus of $\ccM_{g,n}^\ld\ra\ccM_{g,n}$, the induced cover $S\ra C$ is {\'e}tale. Therefore $S$
is a tree of elliptic curves as well. Note that for any circle $\alpha$ on one
of the components of $S$ there is a degeneration of $S$ in the family
$\ccM_{g,n}^\ld\ra\ccM_{g,n-1}^{p(\ld)}$, induced by a degeneration of $C$ in the family
$\ccM_{g,n}\ra\ccM_{g,n-1}$, in which $\alpha$ is a vanishing loop.

From the exact sequence of groups
$\pi_1(S)\ra\pi_1({\ccM}_{g,n}^\ld)\stackrel{\phi}{\ra}\pi_1({\ccM}_{g,n-1}^{p(\ld)})\ra1$, 
given by Lemma~\ref{exact sequence}, it follows that
$\pi_1(\ccM_{g,n}^\ld)=\pi_1(\ccM_{g,n-1}^{p(\ld)})$.
\end{proof}

In the previous section, we defined the Galois level $\GG^{w(k,m)}$ as the kernel of the natural
representation $\rho_{W^{k,m}}:\GG_{g,n}\ra\out(\Pi_{g,n}/W^{k,m}\Pi)$. 

The group $\cK_{g,n}(m)$ is defined, for $2g-2+n>0$, to be the normal subgroup of $\GG_{g,n}$
generated by Dehn twists along separating circles, $m$-th powers of Dehn twists along non-separating
circles and $m_2$-th powers of {\it bounding pair maps} $\tau_{\gm_1}^{-1}\tau_{\gm_2}$, for all cut pairs
$\gm_1,\gm_2$ bounding an unmarked subsurface of $S_{g,n}$ (as usual, we let $m_2:=m/gcd(2,m)$). 

\begin{proposition}\label{generation}Let $g\geq 2$. Then, it holds $\GG_g^{w(3,m)}=\cK_g(m)$.
\end{proposition}
\begin{proof}For $g=2$, the proposition is just Proposition~3.2 in \cite{B-P}. The case $m$ odd was already
treated in Proposition~3.4 in \cite{B-P}. Here, we will briefly recall that proof, pointing out the slight
modifications needed in order to make it work for $m$ even. 

One starts with the remark that $\cK_g(m)\leq\GG_g^{w(3,m)}$ and then observes that
$\cK_g(m)\cdot\cT_g=\GG_g^{w(3,m)}\cdot\cT_g$, since both groups, by Proposition~3.4 in \cite{B-P},
equal $\GG(m)$. Therefore, in order to prove the proposition, it is enough to show that 
$\cK_g(m)\cap\cT_g=\GG_g^{w(3,m)}\cap\cT_g$. Since there is an inclusion
$\cK_g(m)\cap\cT_g\leq\GG_g^{w(3,m)}\cap\cT_g$, it is then enough to show that:
$$(\cK_g(m)\cap\cT_g)/\cK_g=(\GG_g^{w(3,m)}\cap\cT_g)/\cK_g\,\,\,\,\,\,\,\,\,\,\,\,\,(\ast)$$
The advantage of considering the latter identity is that it can be checked inside a torsion free abelian
group. More precisely, there is a natural representation (Johnson's homomorphism):
$$j_0:\cT_g\ra\wedge^3 H_1(S_g,\Z)/([S_g]\wedge H_1(S_g,\Z)),$$
where $[S_g]\in\wedge^{2} H_1(S_g,\Z)$ is the fundamental class of $S_g$. The right-hand side is a
free $\Z$-module of rank ${\st {2g \choose 3}-2g}$ and $\ker j_0=\cK_g$.

From Corollary~6.4 in \cite{P-dJ}, it easily follows that the group $\pi_1(S_g)^{[3],m}$ generates inside the
free abelian group $\pi_1(S_g)^{[2]}/\pi_1(S_g)^{[3]}$ the subgroup of $m_2$-th powers, i.e.
$$\pi_1(S_g)^{[3],m}/\pi_1(S_g)^{[3]}\cong m_2\cdot\wedge^2 H_1(S_g,\Z)/[S_g].$$
Therefore, the image of $H_1(S_g,\Z)\otimes\pi_1(S_g)^{[3],m}$ inside $\wedge^3
H_1(S_g,\Z)/([S_g]\wedge H_1(S_g,\Z))$ equals the submodule $m_2\cdot\wedge^3
H_1(S_g,\Z)/([S_g]\wedge H_1(S_g,\Z))$. Since the group $\cT_g/\cK_g$ embeds as a primitive submodule,
with a basis given by bounding pair maps, in the free $\Z$-module $\wedge^3
H_1(S_g,\Z)/([S_g]\wedge H_1(S_g,\Z))$, it follows that an element of $\cT_g/\cK_g$ belongs to
the submodule of $m_2$-th powers if and only if it can be represented as a product of $m_2$-th powers
of bounding pair maps. This immediately yields the identity $(\ast)$, thus completing the proof of the
proposition.
\end{proof}  

We then have the following generalization of Theorem~3.5 in \cite{B-P}:

\begin{theorem}\label{simply connected}\begin{enumerate}
\item Let $g=1$ and $n\geq 1$. For all $m\geq 2$, there are natural isomorphism
$\pi_1(\ccM_{1,n}^{w(3,m)})\cong\pi_1(\ccM^{(m)}_{1,1})$. Therefore, it holds
$\GG^{w(3,m)}_{1,n}/\cK_{1,n}(m)\cong\pi_1(\ccM^{(m)}_{1,1})$. 
 \item Let $g\geq 2$ and $n\geq 0$. For $m\geq 2$, the level structure $\ccM^{w(3,m)}$ over $\ccM_{g,n}$
is simply connected. Therefore, it holds $\GG^{w(3,m)}=\cK_{g,n}(m)$.
\end{enumerate}
\end{theorem}
\begin{proof}It holds $\out^+(\Pi_{1,1}/W^{3,m}\Pi)=\SL_2(\Z/m)$. Therefore, for $g=1$ and $n=1$,
we have $\GG^{w(3,m)}=\GG(m)\cong\pi_1(\cM^{(m)}_{1,1})$. The statement of the theorem, for $g=1$ and
$n=1$, then follows since the generators of $\cK_{1,1}(m)$ correspond, in the fundamental group of
$\cM^{(m)}_{1,1}$, to small loops around the punctures.

For $g\geq 2$, from Proposition~\ref{generation} and the same arguments of the proof of Proposition~3.3 
in \cite{B-P}, it follows that the level structures $\ccM^{w(3,m)}_{g,n}$ are simply connected for $n=0,1$.

For $n\geq 2$, we proceed by induction on $n$. Let us then assume that the statement of the theorem
has been proved up to $n-1$ and let us prove it for $n$. By Lemma~\ref{comparison lemma} and
Theorem~\ref{monodromy coefficients}, it is enough to prove (with the same notations of the lemma) the
identity:
$$p(\GG^{w(3,m)}_{g,n})=\GG^{w(3,m)}_{g,n-1}.$$ 
The inclusion $p(\GG^{w(3,m)}_{g,n})\leq\GG^{w(3,m)}_{g,n-1}$ is trivial, since the epimorphism
$\Pi_{g,n}\tura\Pi_{g,n-1}$ induces an epimorphism  $\Pi_{g,n}/W^{3,m}\Pi\tura\Pi_{g,n-1}/W^{3,m}\Pi$.

For $g\geq 2$, the reverse inclusion follows, from the fact that, by inductive hypothesis:
$$\GG^{w(3,m)}_{g,n-1}=\cK_{g,n-1}(m)= p(\cK_{g,n}(m))\leq p(\GG^{w(3,m)}_{g,n}).$$

For $g=1$ and all $n\geq 1$, there is a series of natural epimorphism: 
$$\GG^{w(3,m)}_{1,n}/\cK_{1,n}(m)\tura\pi_1(\cM^{(m)}_{1,1})\tura\pi_1(\ccM^{(m)}_{1,1}).$$ 
Since, by inductive hypothesis, it holds
$\GG^{w(3,m)}_{1,n-1}/\cK_{1,n-1}(m)\cong\pi_1(\ccM^{(m)}_{1,1})$
and, moreover, $p(\cK_{1,n}(m))=\cK_{1,n-1}(m)$, it follows that the natural homomorphism 
$p:\GG^{w(3,m)}_{1,n}\ra\GG^{w(3,m)}_{1,n-1}$ is surjective.

In order to complete the proof of the theorem, let us just remark that, as shown in the proof of
Corollary~3.11 in \cite{B-P}, the assertion about the fundamental groups of $\ccM^{w(3,m)}$ implies the
assertions about the generators for the corresponding levels.
\end{proof}

Of course, $\hpi_1(\wM)\cong \widehat{\GG/\cK}$ is naturally isomorphic to
$\ilim_{m>0}\,\GG/\cK(m)$. From Theorem~\ref{simply connected}, it then follows:
 
\begin{theorem}\label{algebraic}Let $g\geq 2$ and $n\geq 0$. The algebraic 
fundamental group of $\wM_{g,n}$ is naturally isomorphic to 
$\lim\limits_{\sr{\textstyle\st\longleftarrow}{\sst m>0}}\GG_{g,n}/\GG^{w(3,m)}$, i.e. to the image of
the representation $\hat{\rho}_{W^3}$.
\end{theorem}

\begin{remark}\label{g=1}{\rm For $g=1$ and $n=1$, from Theorem~\ref{simply connected}, it follows that
$$\lim\limits_{\sr{\textstyle\st\longleftarrow}{\sst m>0}}\GG_{1,1}/\GG^{w(3,m)}\cong\SL_2(\ZZ),$$ 
which is not the profinite completion of $\pi_1(\wM_{1,1})\cong\SL_2(\Z)$. In fact, according to
Theorem~8.8.1 in \cite{R-Z}, the kernel of the natural homomorphism $\wh{\SL_2(\Z)}\ra\SL_2(\ZZ)$ is a
free profinite group of countably infinite rank, which we denote by $\hat{F}_\infty$. }
\end{remark} 

Let $2g-2+n>0$, in general, we have $\hpi_1(\wM_{g,n})\cong\hGG_{g,n}/\ol{\cK}_{g,n}$, where $\hGG_{g,n}$
denotes the profinite completion of the Teichm{\"u}ller group $\GG_{g,n}$ and  $\ol{\cK}_{g,n}$ the closure of
the group ${\cK}_{g,n}$ inside $\hGG_{g,n}$. Of course, $\ol{\cK}_{g,n}\leq\ker\hat{\rho}_{W^3}$, therefore
the representation $\hat{\rho}_{W^3}$ induces a natural homomorphism (which we denote in the same way):
$$\hat{\rho}_{W^3}:\hpi_1(\wM_{g,n})\ra\out(\hPi_{g,n}/W^3\hPi).$$
By Theorem~\ref{algebraic}, $\hat{\rho}_{W^3}$ is injective for $g\geq 2$. For $g=1$, it holds the following result:

\begin{theorem}\label{algebraic genus 1}For $n\geq 1$, the kernel of the natural representation:
$$\hpi_1(\wM_{1,n})\ra\out(\hPi_{1,n}/W^3\hPi)$$
is naturally isomorphic to the congruence kernel $\hat{F}_\infty$.
\end{theorem}
\begin{proof}Let us denote by $\ol{\cK}_{1,n}(m)$ the closure of the group ${\cK}_{1,n}(m)$ in the
profinite group $\hGG_{1,n}$. By Theorem~\ref{simply connected}, there is a natural isomorphism
$\hGG_{1,n}^{w(3,m)}/\ol{\cK}_{1,n}(m)\cong\hpi_1(\ccM_{1,1}^{(m)})$. Let us now observe that, identifying
the fundamental group of the level structure $\cM_{1,1}^{(m)}$ with the level $\GG(m)$, there is a natural 
isomorphism
$$\pi_1(\ccM_{1,1}^{(m)})\cong\GG(m)/<\tau_\gm^m|\, \gm\mbox{ a circle on } S_{1,1}>.$$
Therefore, there is a series of natural isomorphisms:
$$\begin{array}{ll}
\lim\limits_{\sr{\textstyle\st\longleftarrow}{\sst m>0}}\hGG^{w(3,m)}_{1,n}/\ol{\cK}_{1,n}(m) &\cong
\lim\limits_{\sr{\textstyle\st\longleftarrow}{\sst m>0}}\hpi_1(\ccM_{1,1}^{(m)})\cong\\
&\cong\lim\limits_{\sr{\textstyle\st\longleftarrow}{\sst m>0}}
\hGG(m)/\ol{<\tau_\gm^m|\, \gm\mbox{ a circle on } S_{1,1}>}
=\lim\limits_{\sr{\textstyle\st\longleftarrow}{\sst m>0}}\hGG(m)=\hat{F}_\infty.
\end{array}$$
The theorem then follows taking the inverse limit on $m$ of the exact sequences of profinite groups:
$$1\ra\hGG^{w(3,m)}_{1,n}/\ol{\cK}_{1,n}(m)\ra\hpi_1(\wM_{1,n})\ra\out(\hPi_{1,n}/W^{3,m}\hPi).$$
\end{proof}

An easy consequence of Theorem~\ref{simply connected} is also:

\begin{proposition}\label{finite}For $g\geq 2$ and $n\geq 0$, let $\wM^\ld$ be a finite connected {\'e}tale
cover of $\wM_{g,n}$, then its compactification $\ccM^\ld$ over $\ccM_{g,n}$ has finite fundamental group.
\end{proposition}
\begin{proof}Let $m$ be the l.c.m. of the ramification indices of the cover $\ccM^\ld\ra\ccM_{g,n}$ over
$\beta_0(\ccM_{g-1,n+2})$ and let $\ol{X}$ be the universal cover of $\ccM^\ld$ and $X$ the inverse image
of $\cM_{g,n}$ in $\ol{X}$. Then it is clear that $\GG_{g,n}\geq\pi_1(X)\geq\cK_{g,n}(2m)$. Therefore,
$\pi_1(X)$ is a finite index subgroup of $\GG_{g,n}$ and $\ol{X}\ra\ccM^\ld$ a finite cover.
\end{proof}
 
Corollary~3.6 in \cite{B-P} is now stated less naively as it follows:

\begin{corollary}\label{naive}Let $\GG^\ld$ be a level in $\GG_{g,n}$ containing $\cK_{g,n}$:
\begin{enumerate}
\item for $g\geq 2$ and $n\geq 0$, it holds $H^1(\ccM^\ld,\Z)=0$;
\item for $g\geq 3$ and $n=0$, it holds $H^1(\GG^\ld,\Z)=0$. 
\end{enumerate}
\end{corollary}
\begin{proof}The first statement immediately follows from Proposition~\ref{finite}. For the second one,
see the proof of Corollary~3.6 in \cite{B-P}.
\end{proof}

The profinite group $\hPi_{g,n}/W^3\hPi$ is nilpotent. Therefore, there is a natural isomorphism:
$$\hPi_{g,n}/W^3\hPi\cong\prod_{\ell\;\mathrm{prime}}\Pi_{g,n}^{(\ell)}/W^3\Pi^{(\ell)},$$
where, for a given prime $\ell>0$, we denote by $\Pi_{g,n}^{(\ell)}$ the pro-$\ell$ completion of $\Pi_{g,n}$.

Let us define $\pi_1^{(\ell)}(\wM_{g,n})$, for a given prime $\ell>0$ and $g\geq 2$, to be the closure of the
image of $\GG_{g,n}$ inside the virtual pro-$\ell$ group $\out(\Pi_{g,n}^{(\ell)}/W^3\Pi^{(\ell)})$. We then
have, for $g\geq 2$:
$$\hpi_1(\wM_{g,n})\cong\prod_{\ell\;\mathrm{prime}} \pi_1^{(\ell)}(\wM_{g,n}).$$

In Corollary~7.2 of \cite{H-L}, the fundamental groups of $\wM_g$ and $\wM_{g,1}$ are explicitly
described. Let us recall this description. Let us denote by $H$ the first integral homology group
of the compact surface $S_g$, obtained filling in the punctures on $S_{g,n}$, and let then
$\omega\in\wedge^2 H$ be the fundamental class of $S_g$. According to some results by Johnson 
(see \cite{J1}, \cite{J2}, \cite{J3}), for $n=0,1$, the Johnson's group $\cK_{g,n}$ is given by the kernel 
of some naturally defined surjective linear representations:
$$j_0:\cT_g\ra\wedge^3 H/(\omega\wedge H) \,\,\,\,\,\mbox{ and }\,\,\,\,\, j_1:\cT_{g,1}\ra\wedge^3 H,$$
called {\it Johnson's homomorphisms}. Let us remark that, for $n=0$, by Johnson's definition, 
$\ker j_0=\ker\rho_{W^3}$. Therefore, the fundamental groups of $\wM_g$ and $\wM_{g,1}$ fit in the
short exact sequences:
$$\begin{array}{l}
1\ra\wedge^3 H/(\omega\wedge H)\ra\pi_1(\wM_{g})\ra\SP_{2g}(\Z)\ra 1\,\\
\\
1\ra\wedge^3 H\ra\pi_1(\wM_{g,1})\ra\SP_{2g}(\Z)\ra 1.
\end{array}$$

More in general, the fundamental group of $\wM_{g,n+1}$ can be described as an abelian extension of the
fundamental group of $\wM_{g,n}$ for all $n\geq 0$:

\begin{theorem}\label{extension}Let $2g-2+n>0$ and $\ell>0$ a prime. 
There are the following natural short split exact sequences:
$$\begin{array}{l}
1\ra\Z^{2g}\ra\pi_1(\wM_{g,n+1})\ra\pi_1(\wM_{g,n})\ra 1\\
\\
1\ra\ZZ^{2g}\ra\hpi_1(\wM_{g,n+1})\ra\hpi_1(\wM_{g,n})\ra 1\\
\\
1\ra\Z_\ell^{2g}\ra\pi_1^{(\ell)}(\wM_{g,n+1})\ra\pi_1^{(\ell)}(\wM_{g,n})\ra 1.
\end{array}$$
\end{theorem}
\begin{proof}The morphism $\wM_{g,n+1}\ra\wM_{g,n}$ is a proper and flat curve endowed with $n$
tautological sections. Let $[C]\in\wM_{g,n}$ be such that the curve $C$ is smooth.
By Lemma~\ref{exact sequence}, there is an exact sequence:
$$\pi_1(C)\ra\pi_1(\wM_{g,n+1})\ra\pi_1(\wM_{g,n})\ra 1.$$
Let us show that the homomorphism $\pi_1(C)\ra\pi_1(\wM_{g,n+1})$ factors through the abelianization
$H_1(C)$ of $\pi_1(C)$. For $g=1$, this is obvious since $\pi_1(C)\cong H_1(C)$. For $g\geq 2$, let us fix a
homeomorphism of the curve $C$ with the reference genus $g$ compact Riemann surface $S_g$. There
is a natural commutative diagram with exact rows:
$$\begin{array}{ccccccc}
1\ra&\pi_1(S_{g,n})&\ra&\GG_{g,n+1}&\ra&\GG_{g,n}&\ra 1\,\\
&\da&&\da&&\da&\\
&\pi_1(S_g)&\ra&\pi_1(\wM_{g,n+1})&\ra&\pi_1(\wM_{g,n})&\ra 1,
\end{array}
$$
where the vertical arrow $\pi_1(S_{g,n})\ra\pi_1(S_g)$ is the epimorphism induced by the inclusion
$S_{g,n}\subset S_g$. Let $\gm$ be a separating circle on $S_{g,n}$ bounding a disc containing all the
punctures of $S_{g,n}$. Let then $S'$ be the genus $g$ subsurface of $S_{g,n}$ with boundary $\gm$ and
let $\GG(S')$ be the mapping class group of $S'$. Fixing the base-point on $\gm$, the fundamental 
group of $S'$ is identified with a subgroup of $\GG(S')$. Moreover, there is a natural monomorphism
$\GG(S')\hookra\GG_{g,n+1}$ which is compatible with the monomorphism $\pi_1(S')\hookra\pi_1(S_{g,n})$
induced by the inclusion $S'\subset S_{g,n}$. By Theorem~2 in \cite{Birman} and Theorem~6 in \cite{J3}, the
quotient of the Torelli subgroup $\cT(S')$ of $\GG(S')$ by the normal subgroup generated by the twists
along separating circles is abelian. In particular, the image of $\cT(S')$ in
$\pi_1(\wM_{g,n+1})$ is abelian. Since the homomorphism
$\pi_1(S')\ra\pi_1(S_g)$ induced by the inclusion $S'\subset S_g$ is an epimorphism and $\pi_1(S')$
is contained in $\cT(S')$, it follows that the image of $\pi_1(S_g)$ in $\pi_1(\wM_{g,n+1})$ is abelian. 
In order to see that this image is actually the abelianization $H_1(S_g,\Z)$, let us observe what follows.

The relative Jacobian $\cJ_{\sst\wM_{g,n+1}/\wM_{g,n}}$, parameterizing relative Cartier divisors 
of degree $0$, is an abelian variety of  rank $2g$ over $\wM_{g,n}$. Therefore,
there is a short exact sequence:
$$1\ra H_1(C,\Z)\ra\pi_1(\cJ_{\sst\wM_{g,n+1}/\wM_{g,n}})\ra\pi_1(\wM_{g,n})\ra 1.$$
For $n\geq 1$, we can associate to the $n$-th tautological section $s_n:\wM_{g,n}\ra\wM_{g,n+1}$ 
the Abel map $j_{s_n}:\wM_{g,n+1}\dashrightarrow\cJ_{\sst\wM_{g,n+1}/\wM_{g,n}}$. This map is 
defined on the points where the curve $\wM_{g,n+1}\ra\wM_{g,n}$ is smooth. Since the stack
$\wM_{g,n+1}$ is regular, the Abel map induces a homomorphism on fundamental groups
$a_\ast:\pi_1(\wM_{g,n+1})\ra\pi_1(\cJ_{\sst\wM_{g,n+1}/\wM_{g,n}})$ and then a commutative
diagram with exact rows:
$$\begin{array}{ccccccc}
&H_1(C,\Z)&\ra&\pi_1(\wM_{g,n+1})&\ra&\pi_1(\wM_{g,n})&\ra 1\,\\
&\,\,\,\,\da{\st\mathrm{id}}&&\,\,\,\,\da{\st a_\ast}&&\,\,\,\,\da{\st\mathrm{id}}&\\
1\ra& H_1(C,\Z)&\ra&\pi_1(\cJ_{\sst\wM_{g,n+1}/\wM_{g,n}})&\ra&\pi_1(\wM_{g,n})&\ra 1.
\end{array}$$
So, since the bottom line is split exact, the upper line is split exact as well. 

For $n=0$, it is still possible to define a homomorphism $\pi_1(\wM_{g,1})\ra\pi_1(\cJ_{\sst\wM_{g,1}/\wM_g})$, 
which fits in a commutative diagram as above. There is a natural morphism from the universal 
curve $\wM_{g,2}\ra\wM_{g,1}$ to the universal curve $\wM_{g,1}\ra\wM_g$ which contracts the rational components of the fibres. Therefore, the relative Jacobian $\cJ_{\sst\wM_{g,2}/\wM_{g,1}}$ is naturally 
isomorphic to the pull-back, along the natural morphism $\wM_{g,1}\ra\wM_g$, of the relative Jacobian  
$\cJ_{\sst\wM_{g,1}/\wM_{g}}$. In particular, there is a natural rational map 
$j:\wM_{g,2}\dashrightarrow\cJ_{\sst\wM_{g,1}/\wM_g}$, well defined 
on the points where the curve $\wM_{g,2}\ra\wM_{g,1}$ is smooth. The induced homomorphism  $j_{\ast}$
on fundamental groups fits in the commutative diagram:
$$\begin{array}{ccc}
\pi_1(\wM_{g,2})&\tura&\pi_1(\wM_{g,1})\,\\
\,\,\,\da{\sst j_{\ast}}&&\da\\
\pi_1(\cJ_{\sst\wM_{g,1}/\wM_{g}})&\tura&\pi_1(\wM_{g}).
\end{array}$$
The sought for homomorphism is then just the composition $j_{\ast}\circ s_{1\ast}$. This completes the proof of the 
theorem. 

Let us conclude observing that, in particular, for all $n\geq 0$, there is a
natural isomorphism $\pi_1(\wM_{g,n+1})\cong\pi_1(J_{\sst\wM_{g,n+1}/\wM_{g,n}})$.

The proof of the algebraic cases is similar.
\end{proof}

In \cite{M1} (but see also \S 2 in \cite{M2}), Morita extended Johnson's homomorphisms $j_0$ and $j_1$ to
the whole Teichm{\"u}ller group. More precisely, let $\frac{1}{2}\wedge^3 H:=\wedge^3 H\otimes\Z[1/2]$.
There is a natural action of the symplectic group $\SP_{2g}(\Z)$ on $\frac{1}{2}\wedge^3 H$. By means of
this action, let us define semidirect products: 
$$\frac{1}{2}\wedge^3 H\rtimes\SP_{2g}(\Z)\,\,\,\,\mbox{ and }\,\,\,\, \frac{1}{2}\wedge^3
H/(\omega\wedge H)\rtimes\SP_{2g}(\Z),$$  
where $\omega$ denotes the orientation class in $H_2(S_g)$. For $g\geq 2$, Morita extended the Johnson's
homomorphisms to natural linear representations:
$$m_0:\GG_g\ra \frac{1}{2}\wedge^3 H/(\omega\wedge H)\rtimes\SP_{2g}(\Z)\,\,\,\,\mbox{ and }\,\,\,\,
m_1:\GG_{g,1}\ra \frac{1}{2}\wedge^3 H\rtimes\SP_{2g}(\Z),$$
whose images are respectively $\pi_1(\wM_g)$ and $\pi_1(\wM_{g,1})$ and have finite index. In particular,
both $\pi_1(\wM_g)$ and $\pi_1(\wM_{g,1})$ are linear groups. Let us observe now that, in the split short
exact sequence of Theorem~\ref{extension}, for $n\geq 1$, the group $\pi_1(\wM_{g,n})$ acts on $\Z^{2g}$
by means of the natural epimorphism $\pi_1(\wM_{g,n})\tura\SP_{2g}(\Z)$ and the tautological action of
$\SP_{2g}(\Z)$ on $\Z^{2g}$. Therefore, if $\pi_1(\wM_{g,n})$ is linear, the same is true for
$\pi_1(\wM_{g,n+1})$. Similar considerations hold for the profinite groups $\pi_1^{(\ell)}(\wM_{g,n})$.
In conclusion, a simple induction yields:

\begin{proposition}\label{arithmetic}Let $2g-2+n>0$. There is an affine group scheme $\cG_{g,n}$ defined
over $\spec(\Z[1/2])$ such that the group $\pi_1(\wM_{g,n})$ embeds in the group of $\Z[1/2]$-valued
points of $\cG_{g,n}$. Similarly, for $g\geq 2$, the profinite group $\pi_1^{(2)}(\wM_{g,n})$ embeds in the 
group of $\Z_2$-valued points of $\cG_{g,n}$ while, for all primes $\ell>2$, there is a natural isomorphism 
of groups $\pi_1^{(\ell)}(\wM_{g,n})\cong\hom(\spec(\Z_\ell),\cG_{g,n})$.
\end{proposition}

By \cite{B-M-S}, the profinite group $\SP_{2g}(\ZZ)$, for $g\geq 2$, is the profinite completion of
$\SP_{2g}(\Z)$. Let $\hat{H}\cong\ZZ^{2g}$ be the profinite completion of $H$. The algebraic fundamental
groups of $\wM_g$ and $\wM_{g,1}$ are then described by the following exact sequences:
$$\begin{array}{l}1\ra\wedge^3 \hat{H}/(\omega\wedge \hat{H})\ra\hpi_1(\wM_{g,})\ra\SP_{2g}(\ZZ)\ra
1\,\\
\\
1\ra\wedge^3 \hat{H}\ra\hpi_1(\wM_{g,1})\ra\SP_{2g}(\ZZ)\ra 1.
\end{array}$$
In particular, both $\pi_1(\wM_g)$ and $\pi_1(\wM_{g,1})$ are residually finite groups. By 
Theorem~\ref{extension}, the same holds also for $n\geq 2$:

\begin{proposition}\label{residual}Let $2g-2+n>0$. The group $\pi_1(\wM_{g,n})$ is residually finite.
\end{proposition}

We can now prove the key result of the paper:

\begin{theorem}\label{johnson's kernel}Let $2g-2+n>0$. The kernel of 
$\rho_{W^3}:\GG_{g,n}\ra\out(\Pi_{g,n}/W^3\Pi)$ is generated by Dehn twists 
along separating circles, i.e. $\ker\rho_{W^3}=\cK_{g,n}$.
\end{theorem}

\begin{proof}The case of genus $0$ is trivial, since $\ker\rho_{W^3}\geq\cK_{0,n}=\GG_{0,n}$.

On the contrary, the case of genus $1$ is quite tricky.
Let $\cpi_1(\wM_{1,n}):=\hpi_1(\wM_{1,n})/\hat{F}_\infty$, for all $n\geq 1$. By Theorem~\ref{extension}
and Theorem~\ref{algebraic genus 1}, for all $n\geq 1$, there is a short exact sequence:
$$1\ra\ZZ\oplus\ZZ\ra\cpi_1(\wM_{1,n+1})\ra\cpi_1(\wM_{1,n})\ra 1.$$
Since $\cpi_1(\wM_{1,1})\cong\SL_2(\ZZ)$ and the natural morphism $\SL_2(\Z)\ra\SL_2(\ZZ)$ is injective,
the above short exact sequence and a simple induction argument yield that the natural morphism
$\pi_1(\wM_{1,n})\ra\cpi_1(\wM_{1,n})$ is injective for all $n\geq 1$. From 
Theorem~\ref{algebraic genus 1}, it then follows that the natural morphism
$\pi_1(\wM_{1,n})\ra\out(\Pi_{1,n}/W^3\Pi)$ is also injective, which is equivalent to
$\ker\rho_{W^3}=\cK_{1,n}$.

For $g\geq 2$, by Proposition~\ref{residual} and Theorem~\ref{algebraic}, the representation $\rho_{W^3}$
induces a monomorphism $\pi_1(\wM_{g,n})\hookra\out(\hPi_{g,n}/W^3\hPi)$ which factors through the
natural homomorphism $\pi_1(\wM_{g,n})\ra\out(\Pi_{g,n}/W^3\Pi)$, which then is injective, i.e.
$\ker\rho_{W^3}=\cK_{g,n}$.
\end{proof}

\begin{corollary}\label{monodromy char}Let $2g-2+n>0$. The topological fundamental group of $\wM_{g,n}$
is naturally isomorphic to the image of the representation $\rho_{W^3}:\GG_{g,n}\ra\out(\Pi_{g,n}/W^3\Pi)$.
\end{corollary}

It is not hard to extend Johnson's definition of $j_0$ to the $n$-pointed case. Let us consider the
central abelian extension:
$$1\ra W^2\Pi/W^3\Pi\ra\Pi_{g,n}/W^3\Pi\ra H\ra 1.$$
Let $f$ be an automorphism of $\Pi_{g,n}$ inducing the identity on $H$, then the assignment
$\gm\mapsto f(\gm)\cdot\gm^{-1}$, for $\gm\in\Pi_{g,n}$, defines a homomorphism $m_n(f):\Pi_{g,n}\ra
W^2\Pi/W^3\Pi$, which factors through $H$. Let us denote by $\aut_H(\Pi_{g,n})$ the group of
automorphisms of $\Pi_{g,n}$ inducing the identity on $H$. It is easy to check that the assignment
$f\mapsto j_n(f)$ is also a homomorphism. It is then defined an abelian representation:
$$m_n:\aut_H(\Pi_{g,n})\ra\hom(H,W^2\Pi/W^3\Pi).$$
Let us observe that $\inn(\Pi_{g,n})$ is contained in $\aut_H(\Pi_{g,n})$. Let us then define:
$$\hom_{ext}(H,W^2\Pi/W^3\Pi):=\hom(H,W^2\Pi/W^3\Pi)/m_n(\inn(\Pi_{g,n})).$$
It is easy to check that $m_n(\inn(\Pi_{g,n}))$ is a primitive subgroup of
$\hom(H,W^2\Pi/W^3\Pi)$. Hence, the abelian group $\hom_{ext}(H,W^2\Pi/W^3\Pi)$ is torsion free as well.

The homomorphism $m_n$ then induces the abelian representation of the Torelli group:
$$J_n:\cT_{g,n}\ra\hom_{ext}(H,W^2\Pi/W^3\Pi),$$
which we call {\it Johnson's homomorphism} and whose kernel coincides with that of the natural
representation $\rho_{W^3}:\GG_{g,n}\ra\out(\Pi_{g,n}/W^3\Pi)$. From Theorem~\ref{johnson's kernel}, 
it follows that the kernel of $J_n$ is precisely the Johnson's subgroup $\cK_{g,n}$. In particular, it holds
$\ker j_i=\ker J_i$, for $i=0,1$. More importantly, the quotient $\cT_{g,n}/\cK_{g,n}:=I_n$ is a
torsion free abelian group which we are going to determine explicitly. By Theorem~\ref{extension}, there 
is an exact commutative diagram:
$$\begin{array}{ccccccc}
&1&&1&&&\\
&\da&&\da&&&\\
1\ra&\Z^{2g}&\ra&\Z^{2g}&\ra &1&\\
&\da&&\da&&\da&\\
1\ra&I_{n+1}&\ra&\pi_1(\wM_{g,n+1})&\ra&\SP_{2g}(\Z)&\ra 1\\
&\da&&\da&&\,\,\da\wr&\\
1\ra&I_n&\ra&\pi_1(\wM_{g,n})&\ra&\SP_{2g}(\Z)&\ra 1\\
&\da&&\da&&\da&\\
&1&&1&&\,1.&
\end{array}$$
A simple induction on $n$ then yields:

\begin{theorem}\label{shortexact}Let $2g-2+n>0$. The quotient $\cT_{g,n}/\cK_{g,n}$ is a free abelian 
group of rank  ${2g\choose 3} {\st +2g(n-1)}$. Therefore:
\begin{enumerate}
\item For $g=1$ and $n\geq 1$, there is a split short exact sequence:
$$1\ra\Z^{2(n-1)}\ra\pi_1(\wM_{1,n})\ra\SL_2(\Z)\ra 1.$$

\item For $g\geq 2$, there is a short exact sequence:
$$1\ra\Z^{\sst {2g\choose 3}+ 2g(n-1)}\ra\pi_1(\wM_{g,n})\ra\SP_{2g}(\Z)\ra 1.$$
\end{enumerate}
Similar results hold for the profinite groups $\hpi_1(\wM_{g,n})$ and $\pi_1^{(\ell)}(\wM_{g,n})$.
\end{theorem}

Let $2g-2+n>0$ and $g\geq 1$. Consider the following exact commutative diagram (for a given group $G$,
we denote  by $H_1(G)$ its first integral homology group, i.e. its abelianization):
$$\begin{array}{ccccccc}
&1&&1&&1&\\
&\da&&\da&&\da&\\
1\ra&\Pi_{g,n}\cap\cK_{g,n+1}&\ra&\cK_{g,n+1}&\ra &\cK_{g,n}&\ra 1\\
&\da&&\da&&\da&\\
1\ra&\Pi_{g,n}&\ra&\cT_{g,n+1}&\ra&\cT_{g,n}&\ra 1\\
&\da&&\da&&\da&\\
1\ra&H_1(\Pi_{g})&\ra &I_{n+1}&\ra&I_n&\ra 1\\
&\da&&\da&&\da&\\
&1&&1&&\,1.&
\end{array}$$
In the canonical morphism $\Pi_{g,n}\hookra\GG_{g,n+1}$ (obtained, with the notations of
Section~\ref{level}, identifying $\Pi_{g,n}$ with the fundamental group of $S_g\ssm\{P_1,\ldots,P_n\}$ 
based at $P_{n+1}$), a small loop around the puncture $P_i$, for $i=1,\ldots,n$, is sent to the Dehn twist
along a simple circle bounding the $2$-punctured disc containing $P_i$ and $P_{n+1}$ in $S_{g,n}$. Therefore,
all such elements are contained in $\Pi_{g,n}\cap\cK_{g,n+1}$. Thus, passing to abelianizations, we obtain:
$$\begin{array}{ccccccc}
&&&1&&1&\\
&&&\da&&\da&\\
&1&\ra&\ker\pi_{n+1}&\ra&\ker\pi_n&\ra 1\\
&\da&&\da&&\da&\\
1\ra&H_1(\Pi_{g})&\ra& H_1(\cT_{g,n+1})&\ra&H_1(\cT_{g,n})&\ra 1\\
&\,\,\,\da\wr&&\,\,\,\,\,\,\da{\st\pi_{n+1}}&&\,\,\,\,\,\da{\st\pi_n}&\\
1\ra&H_1(\Pi_{g})&\ra &I_{n+1}&\ra&I_n&\ra 1\\
&\da&&\da&&\da&\\
&1&&1&&\,1.&
\end{array}$$ 
Let us observe that, since $I_{n+1}$ and $I_n$ are free abelian groups, all the short exact sequences in the
above commutative diagram are split. Therefore, we have:

\begin{theorem}\label{abelianization}Let $2g-2+n>0$ and $g\geq 1$. For all $n\geq 0$, there is a split natural
short exact sequence:
$$1\ra H_1(\Pi_g,\Z)\ra H_1(\cT_{g,n+1},\Z)\ra H_1(\cT_{g,n},\Z)\ra 1.$$
Therefore, for all $n\geq 0$, the abelianization of the Torelli group $\cT_{g,n}$ is isomorphic to: 
$$H_1(\cT_g,\Z)\oplus H_1(\Pi_g,\Z)^{\oplus n}.$$
\end{theorem}

\begin{remark}\label{abelianization 2}{\rm The above theorem, together with Johnson's computation
of $H_1(\cT_g,\Z)$, for $g\geq 3$, (see \cite{J4}) and Mess' computation of $H_1(\cT_2,\Z)$ (see 
\cite{Mess}) provides a complete description of the abelianization of $\cT_{g,n}$ for all $2g-2+n> 0$.}
\end{remark}

\section[The effect of clutching morphisms]{The effect of the clutching morphisms\\ on fundamental
groups}\label{stratification} 
Let $\sg:=\{\gm_1,\ldots,\gm_k\}$ be a set of distinct, non-trivial, isotopy classes of separating circles on
$S_{g,n}$, such that they admit a set of disjoint representatives none of them bounding a disc with a single
puncture. There is then a homeomorphism $S_{g,n}\ssm\sg\cong\coprod_{i=1}^k S_{g_i,n_i+1}$, where
$\sum g_i=g$, $\sum n_i=n$ and $2g_i-1+n_i>0$ for $i=0,\ldots,k$. Every such set
$\sg$ determines a clutching morphism ({\'e}tale on the image):
$$\dd_\sg:\wM_{g_0,n_0+1}\times\ldots\times \wM_{g_k,n_k+1}\ra\wM_{g,n},$$
which is an embedding except for $n=0$, $k=1$ and $g_0=g_1$. The image of $\dd_\sg$ is the closed
irreducible substack of $\wM_{g,n}$ who generically parametrizes singular curves homeomorphic to
the surface obtained collapsing on $S_{g,n}$ the circles in the set $\sg$. 

The aim of this section is to describe the effect of the above natural morphisms on fundamental groups. 

\begin{theorem}\label{stratum}Let $2g-2+n>0$ and let $\sg$ be a set of circles on $S_{g,n}$ as above. The
clutching morphism $\dd_\sg$ induces on topological fundamental groups a monomorphism:
$$\dd_{\sg\ast}:\pi_1(\wM_{g_0,n_0+1})\times\ldots\times\pi_1(\wM_{g_k,n_k+1})
\hookra\pi_1(\wM_{g,n}).$$
\end{theorem}

\begin{proof}It is clearly enough to prove the theorem for $k=1$. Let then $\gm$ be a separating circle
on $S_{g,n}$ such that there is a homeomorphism $S_{g,n}\cong S_{g_0,n_0+1}\coprod S_{g_1,n_1+1}$.
In order to prove that the induced homomorphism $\dd_{\gm\ast}$ is injective, we need to introduce
some notations.

Let $\hM_{g,n}$ be the real oriented blow-up of $\wM_{g,n}$ along the divisor $\wM_{g,n}\ssm\cM_{g,n}$
(see Section~3 of \cite{B}, for more details on this construction). There is a natural embedding
$\cM_{g,n}\hookra\hM_{g,n}$ which is a homotopy equivalence. So, the natural projection
$\hM_{g,n}\ra\wM_{g,n}$ induces a short exact sequence:
$$1\ra\cK_{g,n}\ra\pi_1(\hM_{g,n})\ra\pi_1(\wM_{g,n})\ra 1.$$
Let $\hd_\gm\ra\hM_{g,n}$ be the pull-back of the clutching morphism
$\cM_{g_0,n_0+1}\times\cM_{g_1,n_1+1}\ra\wM_{g,n}$ along the projection
$\hM_{g,n}\ra\wM_{g,n}$. The morphism $\hd_\gm\ra\hM_{g,n}$ induces a monomorphism
$\pi_1(\hd_\gm)\hookra\GG_{g,n}$ on fundamental groups, which identifies $\pi_1(\hd_\gm)$ with the
stabilizer $\GG_{\vec{\gm}}$ of the oriented isotopy class of the circle $\gm$ on $S_{g,n}$, for the natural
action of the mapping class group $\GG_{g,n}$. The natural projection
$\hd_\gm\tura\cM_{g_0,n_0+1}\times\cM_{g_1,n_1+1}$ is the $S^1$-bundle inducing the classical
central extension which describes the stabilizer $\GG_{\vec{\gm}}\cong\pi_1(\hd_\gm)$:
$$1\ra\Z\cdot\tau_\gm\ra\pi_1(\hd_\gm)\sr{p}{\ra}\GG_{g_0,n_0+1}\times\GG_{g_1,n_1+1}\ra 1.$$
There is also a short exact sequence:
$$1\ra\cK_{g_0,n_0+1}\times\cK_{g_1,n_1+1}\ra\GG_{g_0,n_0+1}\times\GG_{g_1,n_1+1}\ra
\wM_{g_0,n_0+1}\times\wM_{g_1,n_1+1}\ra 1.$$
From the commutative diagram:
$$\begin{array}{ccc}
\GG_{\vec{\gm}}&\hookra &\GG_{g,n}\\
\da^p&&\da\\
\GG_{g_0,n_0+1}\times\GG_{g_1,n_1+1}&\ra&\pi_1(\wM_{g,n}),
\end{array}$$
it then follows that, in order to prove the theorem, all we need to show is that 
$$p(\GG_{\vec{\gm}}\cap\cK_{g,n})=\cK_{g_0,n_0+1}\times\cK_{g_1,n_1+1}.$$
The inclusion $\cK_{g_0,n_0+1}\times\cK_{g_1,n_1+1}\leq p(\GG_{\vec{\gm}}\cap\cK_{g,n})$ being trivial, it
is enough to prove that:
$$p(\GG_{\vec{\gm}}\cap\cK_{g,n})\leq\cK_{g_0,n_0+1}\times\cK_{g_1,n_1+1}.$$
Let $f\in\GG_\gm$, it can be represented by the product of two commuting homeomorphisms $f_0$
and $f_1$, supported respectively on $S_{g_0,n_0+1}$ and $S_{g_1,n_1+1}$. It is easily checked that
$f\in\cT_{g,n}\geq\cK_{g,n}$ if and only if $p(f_i)\in\cT_{g_i,n_i+1}$, for $i=0,1$. In order to prove that 
$f\in\GG_{\vec{\gm}}\cap\cK_{g,n}$ implies $p(f_i)\in\cK_{g_i,n_i+1}$, for $i=0,1$, by
Theorem~\ref{johnson's kernel}, it is then enough to show that, if $f$ is in the kernel of the Johnson's
homomorphism:
$$J_n:\cT_{g,n}\ra\hom_{ext}(H_1(S_g),W^2\Pi/W^3\Pi),$$
then $p(f_i)$ is in the kernel of the Johnson's homomorphism: 
$$J_{n_i+1}:\cT_{g_i,n_i+1}\ra\hom_{ext}(H_1(S_{g_i,n_i+1}),W^2\Pi/W^3\Pi),\,\,\,\,\,\,\mbox{ for } i=0,1.$$

The embedding $S_{g_i,n_i+1}\hookra S_{g,n}$ induces a monomorphism $\Pi_{g_i,n_i+1}\hookra\Pi_{g,n}$,
for $i=0,1$. By Lemma~2.5 in \cite{B-P} (more directly, from its proof), it follows that 
$$\Pi_{g_i,n_i+1}\cap W^3\Pi_{g,n}=W^3\Pi_{g_i,n_i+1},\,\,\,\,\,\,\,\mbox{ for }i=0,1.$$ 
Therefore, for $i=0,1$, there is a monomorphism $\Pi_{g_i,n_i+1}/W^3\Pi\hookra\Pi_{g,n}/W^3\Pi$ and then
a commutative diagram with exact rows:
$$\begin{array}{ccccc}
1\ra &W^2\Pi_{g_i,n_i+1}/W^3\Pi&\ra\Pi_{g_i,n_i+1}/W^3\Pi\ra&H_1(S_{g_i,n_i+1})&\ra 1\,\\
&\da&\da&\da&\\
1\ra &W^2\Pi_{g,n}/W^3\Pi&\ra\Pi_{g,n}/W^3\Pi\ra&H_1(S_{g})&\ra 1,
\end{array}$$
and the vertical arrows are the monomorphisms which are induced by the embedding $S_{g_i,n_i+1}\hookra
S_{g,n}$. With the notations of Section~\ref{compact type}, it is then clear that, for a given $\alpha\in
\Pi_{g,n}$, it holds  $m_n(\inn(\alpha))(H_1(S_{g_i,n_i+1}))\leq W^2\Pi_{g_i,n_i+1}/W^3\Pi$, if and only if
$\alpha$ actually lies inside $\Pi_{g_i,n_i+1}$. Therefore, $f\in\ker J_n$ if and only if 
$f_i\in\ker J_{n_i+1}$, for $i=0,1$.
\end{proof}

The construction of Johnson's homomorphism has a profinite (and actually a pro-$\ell$ as well)
analogue. Let us sketch the profinite case. The profinite group $\hPi_{g,n}/W^3\hPi$ is an extension of
torsion free abelian groups:
$$1\ra W^2\hPi/W^3\hPi\ra\hPi_{g,n}/W^3\hPi\ra \hat{H}\ra 1.$$
Proceeding as above, let us define abelian representations:
$$\hat{m}_n:\aut_{\hat{H}}(\hPi_{g,n})\ra\hom(\hat{H},W^2\hPi/W^3\hPi)\,\mbox{ and then }\,
\hat{J}_n:\hat{\cT}_{g,n}\ra\hom_{ext}(H,W^2\Pi/W^3\Pi).$$
Since it holds $W^2\hPi/W^3\hPi\cong (W^2\Pi/W^3\Pi)\otimes\ZZ$ and $\hat{H}\cong H\otimes\ZZ$, the 
embedding $S_{g_i,n_i+1}\hookra S_{g,n}$ induces a monomorphism
$\hPi_{g_i,n_i+1}/W^3\hPi\hookra\hPi_{g,n}/W^3\hPi$ as well. It also holds:
$$\hom(\hat{H},W^2\hPi/W^3\hPi)\cong\hom(H,W^2\Pi/W^3\Pi)\otimes\ZZ\,\mbox{ and }\,
\hat{m}_n(\inn(\hPi_{g,n}))\cong m_n(\inn(\Pi_{g,n}))\otimes\ZZ.$$
The same argument of the proof of Theorem~\ref{stratum} then yields the following algebraic version:

\begin{theorem}\label{stratum alg}Let $2g-2+n>0$ and let $\sg$ be a set of circles on $S_{g,n}$ as above. 
The clutching morphism $\dd_\sg$ induces on algebraic fundamental groups a homomorphism:
$$\dd_{\sg\ast}:\hpi_1(\wM_{g_0,n_0+1})\times\ldots\times\hpi_1(\wM_{g_k,n_k+1})
\hookra\hpi_1(\wM_{g,n}),$$
which is injective for $g\leq 1$ and for $g_i\neq 1$, for all $i=0,\ldots,k$. If $g_{i_j}=1$ for $0\leq i_j\leq
k$ and $j=1,\ldots,h$, then the kernel of $\dd_{\sg\ast}$ is the product of the subgroups $\hat{F}_\infty$ 
contained in the factors $\hpi_1(\wM_{1,n_{i_j}+1})$, for $j=1,\ldots k$. A similar result holds for the
profinite group $\pi_1^{(\ell)}(\wM_{g,n})$, for all primes $\ell>0$.
\end{theorem}

Let us remark the analogy with the properties of the natural stratification of the
moduli stack $\cA_g$ of principally polarized complex abelian varieties of dimension $g$, given by the loci of
decomposable abelian varieties. For every couple of positive integers $g_1$ and $g_2$ such that
$g_1+g_2=g$, there is indeed a natural morphisms of D--M stacks:
$$\cA_{g_1}\times\cA_{g_2}\ra\cA_g,$$ 
which is {\'e}tale on the image and an embedding if and only if $g_1\neq g_2$. The homomorphism induced on
topological fundamental groups is always a monomorphism (corresponding to the natural monomorphism
$\SP_{2g_1}(\Z)\times\SP_{2g_2}(\Z)\hookra\SP_{2g}(\Z)$). The homomorphism induced on algebraic
fundamental groups is a monomorphism if $g_i\neq 1$, for both $i=1,2$. Otherwise, its kernel is the product
of the kernels of the natural epimorphisms $\widehat{\SP}_{2g_i}(\Z)\tura\SP_{2g_i}(\ZZ)$,
for $i=1,2$, which is non-trivial if and only if $g_i=1$. A significant (and non-casual) discrepancy is given by
the fact that $\pi_1(\cA_g)$ is non-trivial whenever $\dim\,\cA_g>0$ while $\pi_1(\wM_{0,n})=\{1\}$ for all 
$n\geq 3$. This is strictly linked with the fact that the stack of polarized abelian varieties $\cA_g$ is a
classifying space for the group $\SP_{2g}(\Z)$, while the moduli stack of curves of compact type
$\wM_{g,n}$ is not, in general, a classifying space for $\pi_1(\wM_{g,n})$.

\section{On the cohomology of $\pi_1(\wM_{g,n})$}\label{cohomology}
The Teichm{\"u}ller space $T_{g,n}$ is contractible, hence the cohomology of
the group $\GG_{g,n}$ is naturally isomorphic to the cohomology of the space $\cM_{g,n}$.

Let us define some natural cohomology classes on the Deligne-Mumford compactification $\ccM_{g,n}$.
We will then get cohomology classes for $\cM_{g,n}$ (and then $\GG_{g,n}$) by restriction. 

The morphism which forgets the last labeled point $p:\ccM_{g,n+1}\ra\ccM_{g,n}$ is naturally
identified with the universal curve. Let $\omega_p$ be its relative dualizing line bundle, let the $s_i$, for
$i=1,\ldots,n$ be the tautological sections of $p$ and let $L_i:=s_i^\ast(\omega_p)$. 

The first Chern class of the line bundle $L_i$, which is denoted by $\psi_i$, is a cohomology class in
$H^2(\ccM_{g,n},\Q)$, for $i=1,\ldots,n$. Then, for $r\geq 0$, let
$$\kappa_r:=\pi_\ast(\psi_{n+1}^{r+1})\in H^{2r}(\ccM_{g,n},\Q).$$
The $\psi_i$, for $i=1,\ldots,n$, together with the $\kappa_r$, for $r\geq 0$, are called {\it tautological
classes}. Let us denote in the same way their restrictions to $\cM_{g,n}$. They span a subalgebra of the
cohomology ring $H^\ast(\cM_{g,n},\Q)$, called the {\it tautological algebra}. Its importance comes from 
the fact that all geometrically relevant cohomology classes are contained there. Moreover, recently, it has
been proved (see \cite{M-W}) that in the stable range (see \cite{Harer}), i.e. in degrees less or equal than
$(g-1)/2$, the cohomology of $\cM_{g,n}$ coincides with the tautological algebra.

The natural epimorphism $\GG_{g,n}\tura\pi_1(\wM_{g,n})$ induces a homomorphism of cohomology rings
$H^\ast(\pi_1(\wM_{g,n}))\ra H^\ast(\GG_{g,n})$. It is then natural to ask whether the tautological algebra
is in the image of this homomorphism. As we are going to show, from some results of Morita, it follows 
that this is the case for $g\geq 2$ and all $n\geq 0$. This is in sharp contrast with the case of the natural
epimorphism $\GG_{g,n}\ra\SP_{2g}(\Z)$. In fact, the image of the induced homomorphism
$H^\ast(\SP_{2g}(\Z))\ra H^\ast(\GG_{g,n})$ contains only the classes $\kappa_{2i}$, for $i\geq 0$. 

Let us review Morita's results. Denote respectively by $U_g$ and $U_{g,1}$ the abelian groups
$\frac{1}{2}\wedge^3 H/(\omega\wedge H)$ and $\frac{1}{2}\wedge^3 H$ defined in Section~\ref{compact
type}. In \cite{M2}, Morita showed that, for $g\geq 2$ and $n=0,1$, there is a natural homomorphism of
cohomology rings:
$$H^\ast(U_{g,n},\Q)^{\SP}\ra H^\ast(U_{g,n}\rtimes\SP_{2g}(\Z),\Q),$$
which, combined with the homomorphism:
$$H^\ast(U_{g,n}\rtimes\SP_{2g}(\Z),\Q)\ra H^\ast(\GG_{g,n},\Q),$$
induced by Morita's extension of Johnson's homomorphism, provides a homomorphism of cohomology rings:
$$H^\ast(U_{g,n},\Q)^{\SP}\ra H^\ast(\GG_{g,n},\Q),$$
whose image contains the tautological algebra. In particular, since the above homomorphism
factors through the natural homomorphism $H^\ast(\pi_1(\wM_{g,n}),\Q)\ra H^\ast(\GG_{g,n},\Q)$,
it follows that the image of the latter contains as well the tautological algebra, for $n=0,1$.
Morita's result easily extends to the $n$-pointed case, for all $n\geq 0$:

\begin{proposition}\label{morita}Let $g\geq 2$.
\begin{enumerate}
\item The natural homomorphism
$\GG_{g,n}\ra\pi_1(\wM_{g,n})$ induces a homomorphism on cohomology rings:
$$H^\ast(\pi_1(\wM_{g,n}),\Q)\ra H^\ast(\GG_{g,n},\Q),$$
whose image contains the tautological algebra.

\item The natural homomorphism
$\GG_{g,n}\ra\pi_1^{(\ell)}(\wM_{g,n})$ induces a homomorphism of the $\ell$-adic cohomology rings:
$$H^\ast(\pi_1^{(\ell)}(\wM_{g,n}),\Q_\ell)\ra H^\ast(\GG_{g,n},\Q_\ell),$$
whose image contains the tautological algebra.
\end{enumerate}
\end{proposition}

\begin{proof}The first item of the proposition has been proved for $n=0,1$.
For $n>1$, let us consider the natural epimorphisms $p_i:\GG_{g,n}\ra\GG_{g,1}$, corresponding to filling in
the set of punctures $\{P_1,\ldots,P_n\}\ssm\{P_i\}$ on $S_{g,n}$, for $i=1,\ldots, n$. Let then $S_{g,n}^i$ 
be the Riemann surface with boundary obtained from $S_{g,n}$ removing a small open disc $D_i$ containing
$P_i$ and  let us denote by $\GG(S_{g,n}^i)$ the mapping class group of $S_{g,n}^i$, for $i=1,\ldots,n$. There
is a natural epimorphism $\GG(S_{g,n}^i)\tura\GG_{g,n}$ with kernel generated by the Dehn twist along the
circle $\dd D_i$. Then, the class $\psi_i\in H^2(\GG_{g,n},\Q)$, for $i=1,\ldots,n$, corresponds 
to the central extension:
$$1\ra\Z\ra\GG(S_{g,n}^i)\ra\GG_{g,n}\ra 1.$$
It is clear that such extension is just the pull-back along $p_i$ of the extension:
$$1\ra\Z\ra\GG(S_{g,1}^1)\ra\GG_{g,1}\ra 1,$$
i.e. that $\psi_i=p_i^\ast(\psi_1)$, for $i=1,\ldots,n$. For $i=1,\ldots,n$, there is a commutative diagram:
$$\begin{array}{ccc}
H^2(\pi_1(\wM_{g,1}))&\ra &H^2(\GG_{g,1})\,\\
\da&&\da{\sst p_i^*}\\
H^2(\pi_1(\wM_{g,n}))&\ra &H^2(\GG_{g,n}).
\end{array}$$
Since the image of $H^2(\pi_1(\wM_{g,1}))$ in $H^2(\GG_{g,1})$ contains $\psi_1$, it follows that 
$p_i^\ast(\psi_1)=\psi_i$ is in the image of $H^2(\pi_1(\wM_{g,n}))$ in $H^2(\GG_{g,n})$, for $i=1,\ldots,n$.

As to the classes $\kappa_r$, for $r\geq 1$, letting $p:\GG_{g,n+1}\ra\GG_{g,n}$ be the epimorphism
corresponding to filling in the $n+1$-th puncture on $S_{g,n+1}$, it holds the identity 
$p^\ast(\kappa_1)=\kappa_1-\psi_{n+1}$ in $H^2(\GG_{g,n+1})$ and, more in general, $p^\ast(\kappa_r)$
equals the sum of the class $\kappa_r$ in $H^{2r}(\GG_{g,n+1})$  plus a polynomial in lower degree
tautological classes (see \cite{A-C}). Therefore, a simple induction and an argument similar to 
the above yields that $\kappa_r\in H^\ast(\GG_{g,n})$ is in the image of $H^\ast(\pi_1(\wM_{g,n}))$, for all
$r\geq 1$, thus completing the proof of item $i.)$ of the theorem.

In order to prove item $ii.)$, we just need to remark the following.
Let $U_{g,n}^{(\ell)}$, for $n=0,1$, be the pro-$\ell$ completion of the finitely generated abelian group
$U_{g,n}$. The natural homomorphism $U_{g,n}\ra U_{g,n}^{(\ell)}$ then induces an isomorphism on 
$\ell$-adic cohomology $H^\ast(U_{g,n}^{(\ell)},\Q_\ell)\cong H^\ast(U_{g,n},\Q_\ell)$, for $n=0,1$.
The argument produced to prove $i.)$ then extends, almost trivially, to the $\ell$-adic case.
\end{proof}

\begin{remark}\label{kawazumi}{\rm In \cite{K-M}, it is actually proved that the image of the natural
homomorphism $H^\ast(\pi_1(\wM_{g,n}),\Q)\ra H^\ast(\GG_{g,n},\Q)$, for $n=0,1$, in the stable range,
equals the tautological algebra (this was done before the proof of Mumford's conjecture in \cite{M-W}).}
\end{remark}

\end{document}